\documentclass[11pt,a4paper]{article}

\usepackage[english,francais]{babel}
\usepackage{amsmath,amssymb}
\usepackage{hyperref}
\usepackage[margin=1.8cm]{geometry}

\usepackage{graphics}

\newtheorem{e-proposition}[theorem]{Proposition}

\newtheorem{e-definition}[theorem]{Definition\rm}

\newtheorem{theoreme}{Th\'eor\`eme}[section]
\newtheorem{lemme}[theoreme]{Lemme}
\newtheorem{proposition}[theoreme]{Proposition}

\newtheorem{definition}[theoreme]{D\'efinition\rm}
\newtheorem{remarque}{\it Remarque}

\def\og{\leavevmode\raise.3ex\hbox{$\scriptscriptstyle\langle\!\langle$~}}
\def\fg{\leavevmode\raise.3ex\hbox{~$\!\scriptscriptstyle\,\rangle\!\rangle$}}

\def\ZZ{\mathbb{Z}}
\def\QQ{\mathbb{Q}}
\def\RR{\mathbb{R}}
\def\CC{\mathbb{C}}

\def\St{\mathrm{S}_3}
\def\Sq{\mathrm{S}_4}
\def\Torus{\mathrm{T}}
\def\pTorus{\mathring{\mathrm{T}}}

\def\SL{\textrm{SL}}
\def\Aff{\textrm{Aff}}
\def\Aut{\textrm{Aut}}

\def\BBB{\mathcal{B}}
\def\CCC{\mathcal{C}}

\begin{document}

\title{Un contre-exemple \`a la r\'eciproque du crit\`ere de Forni pour la positivit\'e des exposants de Lyapunov du cocycle de Kontsevich-Zorich}
\author{Vincent Delecroix et Carlos Matheus}

\selectlanguage{francais}

\maketitle

\begin{abstract}

\selectlanguage{francais}
Nous introduisons deux surfaces \`a petits carreaux respectivement de degr\'e $8$ dans la strate $\Omega \mathcal{M}_3(2,2)$ et de degr\'e $9$ dans la strate $\Omega \mathcal{M}_4(3,3)$. Dans ces deux exemples les dimensions des sous-espaces isotropes de l'homologie engendr\'es par les circonf\'erences des cylindres dans les directions rationelles sont respectivement $2$ et $3$ (ind\'ependemment de la pointe).  Ainsi, dans chacun de ces exemples, le crit\`ere g\'eom\'etrique de G.~Forni pour la non-uniforme hyperbolicit\'e du cocycle de Kontsevich-Zorich ne s'applique pas. Cependant, en utilisant un crit\`ere alg\'ebrique pour les conditions de ``twisting'' et ``pinching'' de \cite{AV1} et \cite{AV2} (voir \cite{MMY}) nous d\'emontrons que dans les deux cas que ce spectre est simple et ne contient aucun exposant nul. En particulier, la positivit\'e du spectre de Lyapunov pour une mesure r\'eguli\`ere n'implique pas l'existence de surface compl\'etement p\'eriodique dans le support de cette mesure dont le sous-espace isotrope engendr\'e par les circonf\'erences des cylindres est de dimension maximale.
\end{abstract}

\selectlanguage{english}
\begin{abstract}
{\bf A counterexample to the reciprocal of Forni criterion about positivity of Lyapunov exponents of the Kontsevich-Zorich cocycle} \vspace{0.5cm} \\
We introduce two square-tiled surfaces, one with $8$ squares inside $\Omega \mathcal{M}_3(2,2)$, and the other with $9$ squares inside $\Omega \mathcal{M}_4(3,3)$, respectively. In these examples, the dimensions of the isotropic subspaces (in absolute homology) generated by the waist curves of the maximal cylinders in any fixed rational direction are $2$ and $3$ respectively. Hence, a geometrical criterion of G.~Forni for the non-uniform hyperbolicity of Kontsevich-Zorich (KZ) cocycle can not be applied to these examples. Nevertheless, we prove that there are no vanishing exponents and the spectrum is simple for these two square-tiled surfaces. In particular, the non-vanishing of exponents of KZ cocycle for a regular measure doesn't imply that the support of this measure contains a completely periodic surface whose waist curves of maximal cylinders generates a Lagrangian subspace in its absolute homology.
\end{abstract}

\selectlanguage{francais}

\section{Introduction}
Les surfaces arithm\'etiques (ou surfaces \`a petits carreaux) sont des cas particuliers de surfaces de translation. Ces derni\`eres ont \'et\'e introduites pour \'etudier la dynamique des billards rationnels. L'espace des modules $\Omega \mathcal{M}_g$ des surfaces de translation de genre $g$ est un orbifold affine complexe; autremement dit il est localement un quotient fini d'un espace vectoriel complexe de dimension fini et les changements de cartes sont des transformations affines. Il est muni d'un fibr\'e plat r\'eel de dimension $2g$ : le fibr\'e de Hodge (la fibre en une surface est la cohomologie r\'eelle de cette surface). La monodromie de ce fibr\'e est appel\'e le \emph{cocycle de Kontsevich-Zorich} (KZ).

L'espace des modules des surfaces de translation poss\`ede une action de $\SL(2,\RR)$. Cette action et le comportement du cocycle KZ le long des orbites de cette action donnent des renseignements sur la dynamique des flots de translations : mesures invariantes, d\'eviation des sommes de Birkhoff, m\'elange faible, \ldots. Nous renvoyons aux survols~\cite{Z1,Z2} et \cite{F1} pour les d\'etails.

Une surface de translation dont la $\SL(2,\RR)$-orbite est ferm\'ee dans l'espace des modules est appel\'ee \emph{surfaces de Veech} et son orbite une courbe de Teichm\"uller. Dans ce cas, l'orbite est isomorphe \`a un quotient de volume fini $\SL(2,\RR) / \Gamma(S)$ o\`u $\Gamma(S)$ est le \emph{groupe de Veech}, i.e. le stabilisateur de la surface de translation $S$ pour l'action de $\SL(2,\RR)$. Une \emph{surface \`a petits carreaux} est un rev\^etement fini du tore plat. Toute surface \`a petits carreaux est une surface de Veech et son groupe de Veech est commensurable \`a $\SL(2,\ZZ)$. Nous renvoyons au survol~\cite{HS} pour plus de d\'etails.

L'objet de cette note est la description de deux exemples de surfaces arithm\'etiques qui montrent qu'il n'y a pas de r\'eciproque au crit\`ere g\'eom\'etrique de G.~Forni \cite{F2} pour la non-uniforme hyperbolicit\'e du cocycle de Kontsevich-Zorich ($KZ$). 

\begin{theoreme}[\cite{F2}, Theorem~1.8] \label{thm:Forni_criterion}
Soit $\Omega \mathcal{M}_g(\kappa)$ une strate de l'espace des modules des diff\'erentielles ab\'eliennes sur les surfaces de genre $g$. Soit $\mu$ une mesure affine $\SL(2,\RR)$-invariante et ergodique. Si dans le support de $\mu$ il existe une surface compl\'etement p\'eriodique $S$ dont les circonf\'erences des cylindres engendrent un sous-espace isotrope de dimension $k$ de l'homologie alors au moins $k$ des exposants de Lyapunov du cocycle $KZ$ pour la mesure $\mu$ sont positifs.
\end{theoreme}

Soit $\Torus$ le tore $\CC / \ZZ^2$ muni de la forme ab\'elienne $dz$. On note $\pTorus$ la surface $\Torus$ priv\'ee de l'origine et $F_2 = \langle r,u\rangle = \pi_1(\pTorus)$ o\`u $r$ et $u$ d\'esigne les g\'en\'erateurs canoniques horizontal et vertical du groupe fondamental.

Soit $\St \in \Omega\mathcal{M}_3(2,2)$ le rev\^etement de degr\'e $8$ de $\pTorus$ dont les holonomies le long de $r$ et $u$ sont respectivement donn\'ees par (voir aussi la figure~\ref{fig:S3_homology})
$$
\bar{r} = (1,2,3,4)(5,6,7,8)
\quad \mathrm{et} \quad
\bar{u} = (1,2,3,5)(4,8,7,6).
$$
Le commutateur de $\bar{r}$ et $\bar{u}$ vaut $[\bar{r},\bar{u}] = \bar{r} \bar{u} \bar{r}^{-1} \bar{u}^{-1} = (2,4,5)(3,8,6)$ (nous prenons ici la convention de multiplier les permutations de gauche \`a droite). Remarquons que l'action de $M=\langle\bar{r},\bar{u}\rangle$ sur $\{1,2,3,4,5,6,7,8\}$ poss\`ede deux blocs $\{1,3,6,8\}, \{2,4,5,7\}$ et montre ainsi que $\St$ est un rev\^etement de degr\'e $4$ d'un tore ramifi\'e au-dessus de deux points.

Rappelons que dans toutes les directions de pentes rationnelles, une surface \`a petits carreaux se d\'ecompose en un nombre fini de cylindres. C'est par exemple le cas des directions horizontales et verticales. Les classes d'\'equivalences de ces d\'ecompositions sont donn\'ees par les \emph{pointes} (ou "\emph{cusp}" en anglais) de la courbe de Teichm\"uller; autrement dit les classes de conjugaisons d'\'el\'ements paraboliques primitifs dans le groupe de Veech.

Dans cette note, nous d\'emontrons les deux r\'esultats suivants. Premi\`erement, les espaces isotropes associ\'es aux d\'ecompositions en cylindres de la surface $\St$ sont de dimension strictement plus petite que le genre de $\St$ qui est $3$:
\begin{proposition} \label{prop:decompositions_cylindres_S3}
La courbe de Teichm\"uller $\SL(2,\RR) \cdot \St$ poss\`ede deux pointes. Dans chacune des directions de pointes, les circonf\'erences des cylindres engendrent un espace isotrope de dimension $2$ dans $H_1(\St;\ZZ)$.
\end{proposition}

Deuxi\`emement, m\^eme si le crit\`ere de Forni (th\'eor\`eme~\ref{thm:Forni_criterion}) ne permet pas de d\'emontrer la positivit\'e des exposants dans ce cas, nous obtenons un th\'eor\`eme de positivit\'e et de s\'eparation des exposants:
\begin{theoreme} \label{thm:positivite_separation_S3}
Les exposants du cocycle $KZ$ pour la mesure de Haar sur $\SL(2,\RR) \cdot \St$ v\'erifient $1 = \nu_1 > \nu_2 > \nu_3 > 0$.
\end{theoreme}

\begin{remarque} Consid\'erons $\Sq \in \Omega \mathcal{M}_4(3,3)$ le rev\^etement de $\pTorus$ donn\'e par 
$$
\bar{r} = (1,2,3)(4,5,6)(7,8,9)
\quad \mathrm{et} \quad
\bar{u} = (4,2,1)(3,6,9)(7,8,5)
$$
La m\'ethode de cette note permet de prouver que, comme pr\'ec\'edemment, les espaces isotropes associ\'es aux d\'ecompositions en cylindres de la surface $\Sq$ sont de dimensions strictement plus petite que le genre de $\Sq$ (\'egal \`a 4) et, dans ce cas aussi, il y a positivit\'e et s\'eparation des exposants. Plus pr\'ecis\'ement, il est possible de montrer que la courbe de Teichm\"uller $\SL(2,\RR) \cdot \Sq$ poss\`ede $8$ pointes, dans chacune des directions de pointes les circonf\'erences des cylindres engendrent un sous-espace isotrope de dimension $3$ dans $H_1(\Sq; \ZZ)$ et les exposants de Lyapunov du cocycle $KZ$ pour la mesure de Haar sur $\SL(2,\RR) \cdot \Sq$ v\'erifient $1 = \nu_1 > \nu_2 > \nu_3 > \nu_4 > 0$.
\end{remarque}




Les surfaces $\St$ et $\Sq$ ont \'et\'e obtenues par g\'en\'eration exhaustive des orbites de surfaces \`a petits carreaux. Il y a par exemple 14 $\SL(2,\ZZ)$-orbites de surfaces \`a 8 carreaux dans $\mathcal{M}_3^{odd}(2,2)$ et 5 orbites pour 9 carreaux dans $\mathcal{M}^{nonhyp}_4(3,3)$. Toutes les surfaces \`a petits carreaux jusqu'\`a 10 carreaux ont \'et\'es test\'ees. Ces calculs ont \'et\'e r\'ealis\'es avec le logiciel Sage \cite{Sa} dont le code sp\'ecifique aux surfaces \`a petits carreaux a \'et\'e r\'edig\'e par le premier auteur \`a partir de programmes de S.~Leli\`evre \cite{L} ainsi que la librairie d\'evelopp\'ee par K.~Kremer et G.~Schmith\"usen \cite{KS}.

\section{D\'emonstration des r\'esultats pour $\St \in \Omega\mathcal{M}_3^{odd}(2,2)$}
Pour d\'emontrer la proposition~\ref{prop:decompositions_cylindres_S3} et le th\'eor\`eme~\ref{thm:positivite_separation_S3} on utilise un crit\`ere de \cite{MMY} donnant une condition suffisante pour la positivit\'e et la s\'eparation des exposants.

Soit $S$ une surface de translation. On note $\Aff(S)$ son groupe affine (le groupe des diff\'eomorphismes de $S$ qui sont lin\'eaires dans les cartes) et $\Gamma(S)$ son groupe de Veech (la partie lin\'eaire de $\Aff(S)$). De mani\`ere g\'en\'erale, on a une suite exacte $1 \rightarrow \Aut(S) \rightarrow \Aff(S) \rightarrow \Gamma(S) \to 1$. Dans notre cas, la surface $\St$ ne poss\`ede pas d'automorphisme distinct de l'identit\'e et nous identifierons dans la suite les matrices de $\Gamma(\St)$ avec les \'el\'ements du groupe affine.

Nous aurons besoin des g\'en\'erateurs de $\SL(2,\ZZ)$ suivants
$$
L = \left(\begin{array}{ll}1&1\\0&1\end{array}\right)
\quad \mathrm{et} \quad
R = \left(\begin{array}{ll}1&0\\1&1\end{array}\right).
$$
Soit $S$ une surface \`a petits carreaux sans automorphisme non-trivial et $\pi:S\to T$ le rev\^etement fini associ\'e. On note $H_1(S; \QQ) = H^{st}_1(S;\QQ) \oplus H_1^{(0)}(S;\QQ)$ o\`u le sous-module $H_1^{(0)}(S;\QQ)$ est l'ensemble des vecteurs d'holonomie nulle pour la forme $\pi^* dz \in \Omega(S)$ et $H^{st}_1(S;\QQ) \subset H_1(S;\QQ)$ est l'espace engendr\'e par les cycles $\pi^{-1}(\gamma)$ o\`u $\gamma$ est une courbe ferm\'ee sur le tore. Si $\phi \in \Aff(S)$ est un \'el\'ement du groupe affine, alors il pr\'eserve le sous-module $H_1^{(0)}(S;\QQ)$ et on note $\phi^{(0)}$ la restriction de $\phi_*$ \`a $H_1^{(0)}(S;\QQ)$.

\begin{definition}
Soit $S$ une surface \`a petits carreaux sans automorphismes non-triviaux. Un \'el\'ement $\phi \in \Aff(S)$ est appel\'e $b$-r\'eduit si il est de la forme $L^{a_1}R^{a_2}\dots L^{a_{2k-1}}R^{a_{2k}}$ avec $k\geq 1$ et $a_1,\dots, a_{2k}\geq 1$ entiers.
\end{definition}

Cette d\'efinition a \'et\'e introduite par \cite{MMY} motiv\'ee par le codage du flot g\'eod\'esique sur la courbe modulaire $\mathbb{H}/SL(2,\mathbb{Z})$ par l'algorithme de fractions continues (voir la section 3 de \cite{MMY} et, en particulier la d\'efinition 3.7, la proposition 3.8, le corollaire 3.9 et la remarque 3.10 de \cite{MMY}).

On \'enonce maintenant le crit\`ere de \cite{MMY} qui est bas\'e sur les travaux d'A.~Avila et M.~Viana \cite{AV2}. Dans \cite{AV2} est donn\'e un crit\`ere g\'en\'eral de s\'eparation et de positivit\'e des exposants de Lyapunov pour les cocycles symplectiques localement constant. Ce crit\`ere a entre autre permis de d\'emontrer la s\'eparation des exposants de Lyapunov du cocycle $KZ$ pour les mesures de Liouville sur les strates \cite{AV1}.

\begin{theoreme}[\cite{MMY}, Theorem 5.4] \label{thm:MMY_criterion}
Soit $S$ une surface \`a petits carreaux r\'eduite de genre $g$ sans automorphismes non-triviaux. Soit $\phi_1, \phi_2 \in \Aff(S)$ deux \'el\'ements $b$-r\'eduits et $\phi_1^{(0)}, \phi_2^{(0)}: H^{(0)}_1(S;\QQ) \rightarrow H^{(0)}_1(S;\QQ)$ les actions induites sur l'homologie. Soit 
$K_1$, resp. $K_2$, le corps de d\'ecomposition de $\phi_1^{(0)}$, resp. $\phi_2^{(0)}$, dans $\CC$. Si
\begin{enumerate}
  \item le polynome caract\'eristique de $\phi^{(0)}_1$ est irr\'eductible sur $\QQ$, ses racines sont toutes r\'eelles, et son groupe de Galois est isomorphe \`a $S_{g-1}\rtimes \{\pm1\}^{g-1}$, 
\item le polyn\^ome minimal de $\phi^{(0)}_2$ est de degr\'e $>2$ et sans facteur irr\'eductible de degr\'e pair,
\item $K_1\cap K_2=\QQ$,
\end{enumerate}
alors le cocycle $KZ$ au-dessus de la courbe de Teichm\"uller de $S$ a $g$ exposants positifs et s\'epar\'es.
\end{theoreme}
Dans la suite, nous construisons deux \'el\'ements $\phi_1, \phi_2 \in \Gamma(\St)$ v\'erifiant les hypoth\`ese du th\'eor\`eme ci-dessus. Nous commen\c cons par d\'ecrire la courbe de Teichm\"uller de la surface $\St$.
\begin{proposition}\label{p.courbeTeichS3}
Le groupe de Veech $\Gamma(\St)$ de l'origami $\St$ est le sous-groupe d'indice $3$ de $\SL(2,\ZZ)$ donn\'e par
$$
\Gamma(\St) = \left\{\left(\begin{array}{ll}a&b\\c&d\end{array}\right) \in \SL(2,\ZZ);\ a+b \equiv c+d \equiv 1\ \mathrm{mod}\ 2 \right\}.
$$
En particulier, la courbe de Teichm\"uller $\CCC$ est de genre $0$ et poss\`ede deux pointes.
\end{proposition}

\begin{figure}[ht!]
\centering
\includegraphics{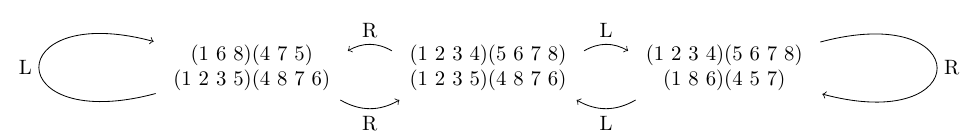}
\caption{Orbite de $\St$ sous l'action de $\SL(2,\ZZ)$.}
\label{fig:S3_orbit}
\end{figure}

\noindent
\textit{D\'emonstration.} Nous appliquons l'algorithme de \cite{Sc1} pour d\'eterminer l'action de $L$ et $R$ sur les classes d'\'equivalences de permutations. Le groupe de Veech est alors donn\'e par l'ensemble obtenu par la composition des matrices $L$ et $R$ le long des lacets enracin\'es en $\St$ dans le graphe orient\'e de la figure~\ref{fig:S3_orbit}.\\

Le groupe $\Gamma(S_3)$ est engendr\'e par $L^2, LRL, R^2, RLR$, ce qui termine la preuve. \hfill $\square$.

\bigskip

\begin{figure}[ht!]
\centering
\includegraphics{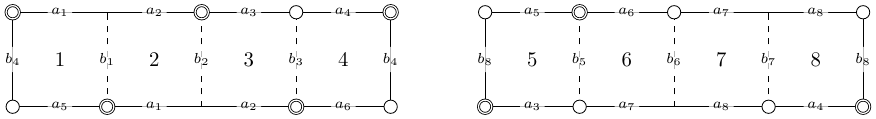}
\caption{La surface $\St$ et sa d\'ecomposition simpliciale.}
\label{fig:S3_homology}
\end{figure}

\begin{figure}[ht!]
\centering
\includegraphics{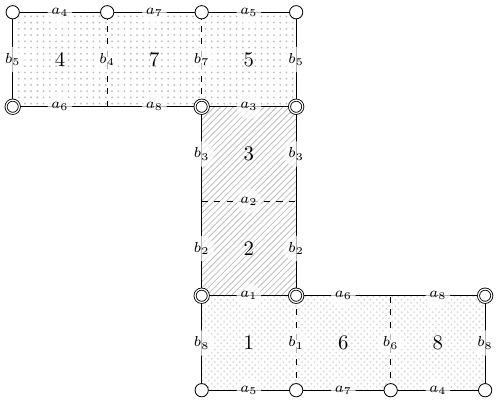}
\caption{La surface $R \cdot \St$ et sa d\'ecomposition simpliciale. Les cylindres form\'es respectivement des carr\'es $\{1,6,8\}$ et $\{4,7,5\}$ sont homologues.}
\label{fig:RS3_homology}
\end{figure}

La proposition \ref{p.courbeTeichS3} permet d'obtenir la proposition \ref{prop:decompositions_cylindres_S3} de la fa\c con suivante:

\noindent
\textit{D\'emonstration de la Proposition \ref{prop:decompositions_cylindres_S3}.} Les d\'ecompositions en cylindres d'une surface compl\'etement periodique dans l'orbite $\SL(2,\RR)\cdot\St$ sont reli\'ees aux pointes de $\Gamma(\St)$, i.e., les $L$-orbites des \'el\'ements de $SL(2,\mathbb{Z})\cdot \St$ (par un lemme d'Anton Zorich, cf. Lemma 2.5 dans \cite{HubLel}). D'apr\`es la proposition \ref{p.courbeTeichS3}, il y a deux d\'ecompositions en cylindres \`a \'etudier reli\'ees aux deux pointes de $\Gamma(\St)$. On v\'erifie (a partir de la figure \ref{fig:S3_homology}) que dans la premi\`ere pointe (i.e., la $L$-orbite de taille $2$ contenant $\St$) la d\'ecomposition est faite de deux cylindres non homologues et (\`a partir de la figure~\ref{fig:RS3_homology}) que dans la seconde (i.e., la $L$-orbite de taille $1$) elle a trois cylindres dont deux sont homologues. \hfill $\square$

\bigskip

Nous choisissons maintenant une base de $H_1(\St;\QQ)$ de mani\`ere \`a pouvoir faire des calculs. La surface $\St$ est munie d'une d\'ecomposition simpliciale canonique $C_1(\St) = \bigoplus_{i=1}^8 \QQ a_i \oplus \bigoplus_{i=1}^8 \QQ b_i$ venant du rev\^etement $\St \rightarrow \Torus$ (voir figure~\ref{fig:S3_homology}).

On remarquera que tous les cycles relatifs $a_i$ (resp. $b_i$) sont d'holonomie $(1,0)$ (resp. $(0,1)$). Une base de $H^{st}_1(\St; \QQ)$ est donn\'ee par les deux vecteurs $a = \sum_{i=1}^8 a_i$ et $b = \sum_{i=1}^8 b_i$. On choisit pour base de $H_1^{(0)}(\St; \QQ)$ les vecteurs suivants $\BBB = \left\{ \frac{1}{2}(a_1 + a_2 - a_7 -a_8),\ (a_4-a_5),\ \frac{1}{2}(b_1+b_2-b_6-b_7),\ (b_4-b_5) \right\}$. Toutes les matrices dans la suite seront exprim\'ees dans la base $\BBB$. En particulier, la forme d'intersection restreinte \`a $H_1^{(0)}(\St; \QQ)$ est donn\'ee par la matrice
$$ \Omega = \left( \begin{array}{rr}
0&M\\
-M&0
\end{array} \right) \mathrm{ avec }\  M = \left(\begin{array}{rr}1&1\\-1&1\\\end{array}\right).
$$

Nous construisons maintenant des \'el\'ements pseudo-Anosov de $\Aff(\St)$ en suivant un chemin dans le graphe associ\'e \`a l'action de $L$ et $R$ sur la $\SL(2,\ZZ)$ orbite de $\St$ (voir figure~\ref{fig:S3_orbit}). On choisit des chemins particuliers qui sont des produits des \'el\'ements paraboliques $L^2$ et $R^2$.

On \'etablit tout d'abord.
\begin{lemme}
Les matrices sur $H_1^{(0)}$ associ\'ees aux transformations paraboliques $L^2$ et $R^2$ sont
$$
\left(\begin{array}{llll}
-1&0&-1&-1\\
0&1&0&0\\
0&0&0&-1\\
0&0&-1&0
\end{array}\right)
\quad \mathrm{et} \quad
\left(\begin{array}{llll}
0&1&0&0\\
1&0&0&0\\
-1&1&-1&0\\
0&0&0&1
\end{array}\right)
$$
\end{lemme}

Soient $\phi_1 = L^8\,R^2\,L^2\,R^2$ et $\phi_2 = L^6\,R^6\,L^4\,R^2$ vus comme des \'el\'ements de $\Aff(\St)$. On v\'erifie alors
\begin{lemme}
Les \'el\'ements $\phi_1$ et $\phi_2$ v\'erifient les hypoth\`eses du th\'eor\`eme \ref{thm:MMY_criterion}.
\end{lemme}
Les polyn\^omes caract\'eristiques de $\phi_1^{(0)}$ et $\phi_2^{(0)}$ sont respectivement $\chi_1 = x^4 - 2 x^3 - 30 x^2 - 2 x + 1$ et $\chi_2 = x^4 + 22x^3 -78 x^2 + 22x + 1$. 
Notons respectivement $K_1$ et $K_2$ les corps de d\'ecomposition de ces polyn\^omes. 
Pour conclure que $K_1 \cap K_2 = \QQ$ on utilise le lemme suivant.
\begin{lemme} \label{lem:corps_quadratiques_intermediaires}
Soit $P = x^4 + a x^3 + b x^2 + a x + 1 \in \ZZ[x]$ un polyn\^ome r\'eciproque irr\'eductible (dans $\ZZ[x]$). Notons $\alpha,\alpha^{-1},\beta,\beta^{-1} \in \CC$ ses quatre racines. Le corps $K=\QQ[\alpha,\beta]$ est le corps de d\'ecomposition de $P$ et notons
\begin{align*}
K'  &=\QQ[\alpha+\alpha^{-1}] = \QQ[\beta + \beta^{-1}] = \QQ[\sqrt{a^2-4(b-2)}] \\
K'' &=\QQ[(\alpha+\beta)(\alpha^{-1}+\beta^{-1})]=\QQ[(\alpha+\beta^{-1}) (\alpha^{-1}+\beta)] = \QQ[\sqrt{(b+2)^2 - 4a^2}] \\
K''' & =\QQ[\alpha \beta^2 + \beta \alpha^{-2} + \alpha^{-1} \beta^{-2} + \beta^{-1} \alpha^2] = \QQ[\sqrt{((b+2)^2-4a^2)(a^2 - 4(b-2))}].
\end{align*}
Alors, $K'$ est de degr\'e $2$. De plus
\begin{enumerate}
\item si $K''=\QQ$, alors $K$ est de degr\'e $4$ de groupe de Galois $V_4=\ZZ/2\ZZ\times\ZZ/2\ZZ$ (groupe de Klein)
\item si $K'''=\QQ$, alors $K$ est de degr\'e $4$ de groupe de Galois $\ZZ/4\ZZ$.
\item si $K'$, $K''$ et $K'''$ sont de degr\'e $2$, alors $K$ est de degr\'e $8$, son groupe de Galois est $D_4$ (groupe dih\'edral) et $K$ poss\`ede exactement trois sous-corps de degr\'e $2$ qui sont $K'$, $K''$ et $K'''$.
\end{enumerate}
\end{lemme}
Ce lemme se v\'erifie 
par un calcul direct pr\'esent\'e dans l'annexe \ref{a.Galois} (voir aussi la section 6 de \cite{MMY}).

Les polyn\^omes $\chi_1$ et $\chi_2$ sont r\'eciproques et irr\'eductibles. Ainsi, gr\^ace au lemme~\ref{lem:corps_quadratiques_intermediaires}, nous savons que les degr\'es de $K_1$ et $K_2$ sont des puissances de $2$. Ainsi, montrer que l'intersection $K_1 \cap K_2$ se r\'eduit \`a $\QQ$ revient \`a montrer qu'il n'ont pas de sous-corps quadratiques communs. Le lemme donne que $K_1$ et $K_2$ sont tous deux de degr\'e $8$ et que leurs sous-corps quadratiques sont
\begin{enumerate}
\item $K'_1 = \QQ[\sqrt{3 \cdot 11}]$, $K''_1 = \QQ[\sqrt{3}]$ et $K'''_1 = \QQ[\sqrt{11}]$,
\item $K'_2 = \QQ[\sqrt{3\cdot 67}]$, $K''_2 = \QQ[\sqrt{3\cdot 5}]$ et $K'''_2 = \QQ[\sqrt{5 \cdot 67}]$.
\end{enumerate}
Ainsi, $K_1$ et $K_2$ n'ont pas de sous-corps quadratiques communs et les conditions du crit\`ere du th\'eor\`eme~\ref{thm:MMY_criterion} sont v\'erifi\'ees. Ce qui termine la preuve du th\'eor\`eme \ref{thm:positivite_separation_S3}.

\bigskip

Les valeurs exp\'erimentales pour les exposants de Lyapunov de $\St$ (obtenues gr\^ace \`a un programme d'A. Zorich) sont $\nu_2 \simeq 0.5879$ et $\nu_3 \simeq 0.0787$. Le fait que la valeur de $\nu_3$ soit petite devant $\nu_2$ explique peut-\^etre pourquoi nous avons \'et\'e oblig\'e de prendre des produits relativement grand pour $\phi_1$ et $\phi_2$.

\bigskip

Remarquons d'autre part que dans les deux exemples pr\'esent\'es dans cette note, les dimensions homologiques sont $g-1$ si $g$ d\'esigne le genre de la surface.

\section*{Remerciements} Les auteurs remercient le rapporteur de l'article pour la relecture attentive d'une version pr\'eliminaire de ce texte.

\appendix

\section{D\'emonstration du lemme \ref{lem:corps_quadratiques_intermediaires}}\label{a.Galois}
Le groupe de Galois d'un polyn\^ome quartique irr\'eductible est un sous-groupe transitif de $S_4$. La liste de sous-groupes \emph{transitifs} de $S_4$ est: 
\begin{itemize}
\item[(a)] $S_4$, $A_4$
\item[(b)] $\langle(1234),(13)\rangle$, $\langle(1324),(12)\rangle$, $\langle(1243),(14)\rangle$
\item[(c)] $\langle(1234)\rangle$, $\langle(1324)\rangle$, $\langle(1243)\rangle$
\item[(d)] $\{\textrm{id},(12)(34),(13)(24),(14)(23)\}$
\end{itemize}
Ici, les $3$ sous-groupes list\'es dans (b) sont conjugu\'es et isomorphes au groupe dih\'edral $D_4$, les  $3$ sous-groupes list\'es dans (c) sont conjugu\'es et isomorphes au groupe cyclique $\mathbb{Z}/4\mathbb{Z}$, et le sous-groupe dans (d) est isomorphe au groupe de Klein $V_4=\mathbb{Z}/2\mathbb{Z}\times\mathbb{Z}/2\mathbb{Z}$.

Dans notre cas, le polyn\^ome quartique irr\'eductible $P(x)=x^4+ax^3+bx^2+ax+1$ est r\'eciproque, de fa\c con que ses racines $\alpha, \alpha^{-1}$, $\beta,\beta^{-1}$ sont calcul\'ees par le changement de variables $u=x+1/x$. En effet, il suffit de remarquer que $P(x)=0$ implique que $u^2+au+(b-2)=0$ et de r\'esoudre l'\'equation $x^2-ux+1=0$. Par cette m\'ethode classique, nous trouvons que
\begin{equation} \label{eq:alphabeta}
\alpha^{\pm1}=\frac{u_+\pm\sqrt{u_+^2-4}}{2} \quad \textrm{and} \quad \beta^{\pm1}=\frac{u_-\pm\sqrt{u_-^2-4}}{2}
\end{equation}
o\`u $u_{\pm}=(-a\pm\sqrt{a^2-4(b-2)})/2 := (-a\pm\sqrt{\Delta_1})/2$. En particulier, le corps de d\'ecomposition $K_P=\mathbb{Q}[\alpha,\beta]$ est de degr\'e 
$4$ ou $8$, et le groupe de Galois de $P$ est $V_4,\mathbb{Z}/4\mathbb{Z}$ ou $D_4$.

Pour la suite, notons que si $\delta_{\pm}:=u_{\pm}^2-4$, alors 
\begin{itemize}
\item $\delta_+\delta_-=(b+2)^2-4a^2:=\Delta_2,$
\item $4(\alpha+\beta)(\alpha^{-1}+\beta^{-1}) = 8 + 2(b-2) - 2\sqrt{\delta_+\delta_-} = 8 + 2(b-2) - 2\sqrt{\Delta_2}$ 
\item $4(\alpha+\beta^{-1})(\alpha^{-1}+\beta) = 8 + 2(b-2) + 2\sqrt{\delta_+\delta_-} = 8 + 2(b-2) + 2\sqrt{\Delta_2} $
\item $8(\alpha \beta^2 + \beta \alpha^{-2} + \alpha^{-1} \beta^{-2} + \beta^{-1} \alpha^2)=-4a(b-2)+8a-4\sqrt{\Delta_1\cdot \Delta_2}$
\end{itemize}
\par
Donc, 
$$K_{P}''=\mathbb{Q}[(\alpha+\beta)(\alpha^{-1}+\beta^{-1})]=\mathbb{Q}[(\alpha+\beta^{-1})(\alpha^{-1}+\beta)]=\mathbb{Q}[\sqrt{\Delta_2}]$$ 
et 
$$K_{P}'''=\mathbb{Q}[\alpha \beta^2 + \beta \alpha^{-2} + \alpha^{-1} \beta^{-2} + \beta^{-1} \alpha^2]=\mathbb{Q}[\sqrt{\Delta_1\cdot \Delta_2}]$$

Pour pouvoir distinguer entre les $3$ cas possibles $V_4,\mathbb{Z}/4\mathbb{Z}$ ou $D_4$ pour le groupe de Galois de $P$, notons que l'expression  
$(\alpha+\beta)(\alpha^{-1}+\beta^{-1})$ est $V_4$-invariante mais pas $\mathbb{Z}/4\mathbb{Z}$-invariante ou $D_4$-invariante. 
Donc, le groupe de Galois de $P$ est isomorphe au groupe de Klein $V_4$ si et seulement si $K_{P}''=\mathbb{Q}$. Dans cette situation, le corps de d\'ecomposition $K=\mathbb{Q}[\alpha,\beta]$ est une extension de degr\'e $4$ de $\mathbb{Q}$. Ceci montre la premi\`ere partie du lemme.

Maintenant, si $K_{P}''\neq\mathbb{Q}$, i.e., $\sqrt{\Delta_2}\notin\mathbb{Q}$, le groupe de Galois est $\mathbb{Z}/4\mathbb{Z}$ ou $D_4$, et nous pouvons distinguer ces $2$ cas car l'expression $\alpha \beta^2 + \beta \alpha^{-2} + \alpha^{-1} \beta^{-2} + \beta^{-1} \alpha^2$ 
est $\mathbb{Z}/4\mathbb{Z}$-invariante mais pas $D_4$-invariante. Donc, le groupe de Galois group de $P$ est $\mathbb{Z}/4\mathbb{Z}$ si et seulement si $K_{P}'''=\mathbb{Q}$ et $K_{P}''\neq\mathbb{Q}$. Ainsi, le corps de d\'ecomposition $K=\mathbb{Q}[\alpha,\beta]$ est aussi une extension de degr\'e $4$ de $\mathbb{Q}$. Ceci montre la deuxi\`eme partie du lemme.

Finalement, si $K_{P}'', K_{P}'''\neq\mathbb{Q}$, i.e., $K_{P}''=\mathbb{Q}[\sqrt{\Delta_2}]$ et $K_{P}'''=\mathbb{Q}[\sqrt{\Delta_1\cdot \Delta_2}]$ sont des extensions quadratiques de $\mathbb{Q}$, alors le groupe de Galois de $P$ est $D_4$ et le corps de d\'ecomposition $K=\mathbb{Q}[\alpha,\beta]$ est une extension de degr\'e $8$ de $\mathbb{Q}$. De plus, les sous-corps de $K$ de degr\'e $2$ sur $\mathbb{Q}$ correspondent aux sous-groupes de $D_4$ d'indice $2$, i.e., d'ordre $4$: il y en a $3$ et il n'est pas difficile de verifier que les sous-corps associ\'es sont $K_{P}'$, $K_{P}''$ et $K_{P}'''$. La preuve du lemme est compl\`ete.


\end{document}